\documentclass[12pt]{amsart}
\usepackage[a4paper,asymmetric]{geometry}
\usepackage{mathscinet}
\usepackage{latexsym}
\usepackage{amsthm}
\usepackage{amssymb}
\usepackage{amsfonts}
\usepackage{amsmath}
\newtheorem{theorem}{Theorem}[section]
\newtheorem{thm}[theorem]{Theorem}

\newtheorem{lem}[theorem]{Lemma}
\newtheorem{proposition}[theorem]{Proposition}

\newtheorem{corollary}[theorem]{Corollary}
\newtheorem{assumption}[theorem]{Assumption}
\theoremstyle{definition}

\newtheorem{defn}[theorem]{Definition}
\theoremstyle{remark}

\newtheorem{rem}[theorem]{Remark}
\numberwithin{equation}{section}

 \DeclareMathAlphabet{\mathpzc}{OT1}{pzc}{m}{it}
 
 \newcommand{\E}{\mathbb{E}}            
 \newcommand{\T}{\mathbb{T}}
 \newcommand{\e}{\varepsilon}
 \newcommand{\p}{\partial}

 \newcommand{\Ll}{\langle}
 \newcommand{\Rr}{\rangle}
 \newcommand{\ph}{\mathpzc{h}}
 \newcommand{\pl}{\mathpzc{l}}
 \newcommand{\N}{\mathbb{N}}
 \newcommand{\R}{\mathbb{R}}
 \newcommand{\Z}{\mathbb{Z}}
 
 \newcommand{\mcl}{\mathcal}
 
 \newcommand{\Be}{\begin{equation}}
 \newcommand{\Ee}{\end{equation}}
 \newcommand{\Bs}{\begin{split}}
 \newcommand{\Es}{\end{split}}
  \newcommand{\Bes}{\begin{equation*}}
 \newcommand{\Ees}{\end{equation*}}
 \newcommand{\BT}{\begin{thm}}
 \newcommand{\ET}{\end{thm}}
 \newcommand{\Bp}{\begin{proof}}
 \newcommand{\Ep}{\end{proof}}
 \newcommand{\BL}{\begin{lem}}
 \newcommand{\EL}{\end{lem}}
 \newcommand{\BP}{\begin{proposition}}
 \newcommand{\EP}{\end{proposition}}
 \newcommand{\BC}{\begin{corollary}}
 \newcommand{\EC}{\end{corollary}}
 \newcommand{\BR}{\begin{rem}}
 \newcommand{\ER}{\end{rem}}
 \newcommand{\BD}{\begin{defn}}
 \newcommand{\ED}{\end{defn}}
 \newcommand{\BI}{\begin{itemize}}
 \newcommand{\EI}{\end{itemize}}
 
 \newcommand{\lbr}{\left(}
 \newcommand{\rbr}{\right)}
 \newcommand{\Lbr}{\left[}
 \newcommand{\Rbr}{\right]}
 \newcommand{\LBr}{\left\{}
 \newcommand{\RBr}{\right\}}

\begin{document}
\title
[Modified log-Harnack inequality and asymptotic strong Feller] {A Modified log-Harnack inequality and asymptotically strong Feller property}
\author[L. Xu]{Lihu Xu}
\address{TU Berlin, Fakult\"{a}t II, Institut f\"{u}r Mathematik,
Str$\alpha \beta$e des 17. Juni 136, D-10623 Berlin, Germany}
\email{xu@math.tu-berlin.de}
\thanks{Partial support by the European Research Council under the European Union's
Seventh Framework Programme (FP7/2007-2013) / ERC grant agreement nr. 258237
is gratefully acknowledged.}
\subjclass[2000]{}
\date{}
\maketitle
\begin{abstract} \label{abstract}
We introduce a new functional inequality, which is a modification of log-Harnack inequality established in \cite{RoWa10} and \cite{W10}, and prove that it implies the asymptotically strong Feller property (ASF). This inequality seems to generalize the criterion for ASF in \cite[Proposition 3.12]{HM06}. As a example, we show by an asymptotic coupling that 2D stochastic Navier-Stokes equation driven by highly degenerate
but \emph{essentially elliptic} noises satisfies our modified log-Harnack inequality.  
\end{abstract} \noindent


 \ \\

 \section{Introduction}
  Dimensional free Harnack inequality was introduced by Wang in \cite{Wan97} to study the diffusions on Riemannian manifolds (see also \cite{ATW06, ATW09} for further development).
 Wang type Harnack inequality has been applied to many research problems such as studying ultracontractivity and functional inequalities (\cite{RoWa03, RoWa03-1, Wan99, Wan01}), short-time behaviors
of infinite-dimensional diffusions (\cite{AiKa01,AiZh02,Ka05}), heat kernel estimates (\cite{BoGeLe01,GoWa01}) and so on.
 In recent years, this inequality has also been established and applied intensively in the study of SPDEs (see e.g. \cite{RoWa03, Wan07, LiWa08, DRW09, ES, W10Ann, WX10} and references within). Let $(P_t)_{t \ge 0}$
  be a Markov semigroup on a Polish space $\mcl X$, this type of Harnack inequality can be formulated as
\Be \label{e:Wang}
 (P_tf)^\alpha (x)\le (P_t f^\alpha)(y) \exp[C_\alpha (t,x,y)],\ \ f\ge 0,
\Ee
 where $\alpha >1$ is a constant, $C_\alpha$ is a positive function on $(0,\infty)\times \mcl X^2$ with $C_\alpha(t, x,x)=0$, which is determined by the underlying stochastic equation.

 On the other hand, in some cases Wang type Harnack inequality is not available, so that the following weaker version (i.e. the log-Harnack inequality)
 \Be \label{LH}   P_t\log f(x)\le \log P_t f(y)+ C(t,x,y), \ \ f\ge 1 \Ee
 becomes an alternative tool in the study. In general, according to \cite[Section 2]{W10}, (\ref{LH}) is the limit version of (\ref{e:Wang}) as $\alpha \to\infty.$
This inequality has been established in \cite{RoWa10} and \cite{W10}, respectively,  for   semi-linear SPDEs with multiplicative noise and the Neumann semigroup on non-convex manifolds. \cite{WWX10} shows that stochastic Burgers equation driven by additive noises satisfies a log-Harnack inequality. As for the research on Wang type Harnack and log-Harnack inequalities on stochastic Navier-Stokes type equations, we refer to \cite{WX10} and \cite{WWX10}.

To our knowledge, nearly all the stochastic systems in the above literatures is forced by \emph{nondegenerate}
noises. It is natural to ask whether a stochastic system with degenerate noises satisfies Wang type
Harnack or log-Harnack inequality. This seems still an open problem.
\ \\

The aim of this paper is to study the case of the systems
driven by highly degenerate noises. For highly degenerate systems, one is usually
not able to prove the strong Feller property (see Example 3.15 of \cite{HM06}).
Since Harnack and log-Harnack inequality
implies strong Feller property (\cite{Wan07}, \cite{RoWa10}), there is no hope to prove these inequalities for
highly degenerate systems. On the other hand,
many dissipative systems such as 2D Navier-Stokes and reaction-diffusion equations driven by highly
degenerate noises (\cite{HM06, HM08}) have asymptotically strong Feller property,
it is natural to ask whether we can establish a functional inequalities which implies the asymptotically
strong Feller property. This is the main motivation that we introduce the modified log-Harnack inequality, which seems to give
a criterion for asymptotically strong Feller property more general than that in \cite[Proposition 3.12]{HM06}. This inequality
also gives some pointwise information on Markov semigroups.

\BD [Modified log-Harnack inequality] \label{d:ModLogHar}
Let $\{P_t\}_{t\geq 0}$ be a Markov semigroup on
a Polish space $\mcl X$, it satisfies a modified log-Harnack inequality if there exist some constants $\alpha>0$, $\beta \geq 0$, $C=C(|x|,|y|)>0$, $\tilde C=\tilde C(|x|,|y|)>0$ and a function $\delta(t) \geq 0$ with $\lim_{t\rightarrow \infty} \delta(t)=0$ such that
\Be \label{e:ModLogHar}
P_t \log f(y) \leq \log P_t f(x)+C |x-y|^{\alpha}+\delta(t) \tilde C |x-y|^\beta ||D \log f||_\infty
\Ee
for any bounded differentiable function $f \geq 1$ and $x, y \in \mcl X$. Moreover, $C$ and $\tilde C$ are both continuous w.r.t. $|x|$ and $|y|$.
\ED
The main results of this paper are the following Theorem \ref{t:ModLogHarASF}, Corollary \ref{p:ModLogHarHaiMat} and Theorem \ref{t:SNSModLogHar}. The first two will be proven in section \ref{s:MLHIVsASF}, while the last one will be shown in section \ref{SExa}.

\BT \label{t:ModLogHarASF}
If a Markov semigroup satisfies a modified log-Harnack inequality, then it is asymptotically strong Feller.
\ET
Hairer and Mattingly gave a criterion for asymptotically strong Feller property as the following.
\BP [Proposition 3.12 in \cite{HM06}] Let $\{P_t\}_{t \geq 0}$ be a Markov semigroup on $\mcl X$. If $\{P_t\}_{t \geq 0}$ satisfies the following inequality
\Be \label{e:HaiMatIne}
|D P_t f(x)| \leq C\left(||f||_\infty+\delta(t) ||Df||_\infty\right)
\Ee
for any bounded differentiable function $f:
\mcl X \rightarrow \R$,
where $C=C(|x|)>0$ and $\delta(t) \geq 0$ with $\lim \limits_{t \rightarrow \infty} \delta(t)=0$, then
$\{P_t\}_{t \geq 0}$ is asymptotically strong Feller.
\EP
The next corollary claims that \eqref{e:ModLogHar} with $\alpha=2$ and $\beta=1$
implies a gradient estimate similar to \eqref{e:HaiMatIne}. Therefore, the modified log-Harnack inequality, in some sense, seems to give a more general criterion for asymptotically strong Feller property. Moreover, \eqref{e:ModLogHar} also provides some pointwise information of the semigroups.

\begin{corollary} \label{p:ModLogHarHaiMat}
If a Markov semigroup $\{P_t\}_{t\geq 0}$ satisfies \eqref{e:ModLogHar} with $\alpha=2$ and $\beta=1$, then
\Be \label{e:HaiMatIne1}
|DP_tf(x)| \leq ||f||^2_\infty+C+\delta(t) \tilde C||D f||_\infty
\Ee
for any bounded differentiable function $f$, where $C=C(|x|)$ and $\tilde C=\tilde C(|x|)$ and $\lim_{t \rightarrow \infty} \delta(t)=0$.
\end{corollary}

The following theorem claims that 2D stochastic Navier-Stokes equation forced by degenerate noises satisfies our modified log-Harnack inequality. We need to emphasize that those degenerate noises have \emph{essential ellipticity effect}, see section 4.5 in \cite{HM06} for more details.

\BT \label{t:SNSModLogHar}
Let $\{P_t\}_{t \geq 0}$ be the Markov semigroup generated by Eq. \eqref{e:XEqn} in section 3 below, a 2D stochastic Navier-Stokes system forced by highly
degenerate noises, then $\{P_t\}_{t \geq 0}$ satisfies a modified log-Harnack inequality, which has the exact form as in Theorem \ref{t:ModLogHI} below. Moreover, $\{P_t\}_{t \geq 0}$ satisfies \eqref{e:HaiMatIne1}.
\ET

{\bf Acknowledgements:} We would like to gratefully thank Feng-Yu Wang for the stimulating and instructive discussions and carefully reading the first draft. We also would like to thank Jiang-Lun Wu for the encouragements. Part of the work was done during visiting Mathematics department of Swansea University.
\ \\

\section{Modified log-Harnack inequality and Asymptotically strong Feller property}
\label{s:MLHIVsASF}
Let us first recall the interesting conception of
\emph{asymptotically strong Feller} property, which was
introduced by Hairer and Mattingly in \cite{HM06}. For the more
details, we refer to \cite{HM06}.
\BD Let $\{d_n\}_n$ be an increasing sequence of pseudo metrics (pp 7. \cite{HM06}) on
a Polish space $\mcl X$. If $\lim_{n \rightarrow \infty} d_n(x,y)=1$ for all $x \neq y$, then
$\{d_n\}$ is called a totally separating system of pseudo metrics for $\mcl X$.
\ED
Given a pseudo metric $d$, for any $d$-Lipschitz continuous function $f: \mcl X \rightarrow \R$,
we define the following semi-norm for $f$:
$$||f||_d=\sup_{\stackrel{x, y\in \mcl X}{x \neq y}} \frac{|f(x)-f(y)|}{d(x,y)}.$$
Given $\mu_1$ and $\mu_2$, two positive finite Borel measures on $\mcl X$ with
\emph{equal} mass, we denote by $\mcl C(\mu_1,\mu_2)$ the set of positive measures on $\mcl X^2$
with marginals $\mu_1$ and $\mu_2$, and define
\Bes
||\mu_1-\mu_2||_d=\inf_{\mu \in \mcl C(\mu_1,\mu_2)} \int_{\mcl X^2} d(x,y) \mu(dx,dy).
\Ees
\BD
[Asymptotically strong Feller property] \label{d:ASF}
A Markov transition semigroup $P_t$ on
a Polish space $\mcl X$ is asymptotically strong Feller at $x$ if there exists a
totally separating system of pseudo metrics $\{d_n\}_n$ for $\mcl X$
and a sequence $t_n>0$ such that
\Be \label{e:ASFDef1}
\inf_{U \in \mcl U_x} \lim \limits_{n \rightarrow \infty} \sup
\sup_{y \in U} ||P_{t_n}(x,\cdot)-P_{t_n}(y, \cdot)||_{d_n}=0,
\Ee
where $\mcl U_x=\{U: \  U \ {\rm is \ the \ neighbourhood \ of \ } x\}$.
We call that $P_t$ satisfies asymptotic strong Feller property if it is asymptotic
strong Feller at each $x \in \mcl X$.
\ED
\BR
If $\mcl X$ has a metric, then the definition \eqref{e:ASFDef1} is equivalent to
\Be \label{e:ASFDef2}
\lim_{r \rightarrow 0} \lim_{n \rightarrow \infty} \sup \sup_{y \in B(x,r)} ||P_{t_n}(x,\cdot)-P_{t_n}(y,\cdot)||_{d_n}=0,
\Ee
where $B(x,r)$ is the ball in $\mcl X$ centered at $x$ with radius $r$ under this metric.
\ER
The following two lemmas are Lemma 3.3 and Corollary 3.5 of \cite{HM06} respectively. The first one is a not difficult consequence of Monge-Kantorovich duality (\cite{Vil03}).
\BL \label{l:DuaD}
Let $d$ be a continuous pseudo metric on a Polish space $\mcl X$ and let
$\mu_1$, $\mu_2$ be two positive measures on $\mcl X$ with \emph{equal} mass. Then
we have
\Bes
||\mu_1-\mu_2||_d=\sup_{{||\varphi||_d=1}} \int_{\mcl X} \varphi(x) (\mu_1-\mu_2)(dx).
\Ees
\EL
\BR \label{r:DuaMu1-Mu2}
By \cite{Vil03} (pp. 34), if $d$ is bounded we have
 \Bes
||\mu_1-\mu_2||_d=\sup_{\stackrel{||\varphi||_d=1}{||\varphi||_\infty \leq ||d||_\infty}} \int_{\mcl X} \varphi(x) (\mu_1-\mu_2)(dx).
\Ees
\ER
\BL
Let $\mcl X$ be a Polish space and let $\{d_n\}$ be a totally separating
system of pseudo metrics for $\mcl X$. Then, $||\mu_1-\mu_2||_{TV}=\lim_{n \rightarrow \infty} ||\mu_1-\mu_2||_{d_n}$
for any two positive measures $\mu_1$ and $\mu_2$ with equal mass on $\mcl X$.
\EL
The following theorem is due to Hairer and Mattingly \cite{HM06}.
\BT
If $P_t$ is an asymptotically strong Feller Markov semigroup
and there exists a point $x$ such that $x \in supp(\mu)$ for every invariant
probability measure $\mu$ of $P_t$, then there exists at most one invariant probability
measure for $P_t$.
\ET

Hairer and Mattingly proved by the above theorem, together with \eqref{e:HaiMatIne}, the following important result: For 2D stochastic Navier-Stokes systems if at least two linearly independent
Fourier modes with different Euclidean norms are driven, then the system is ergodic (see \cite[Theorem 2.1]{HM06} and the examples therein). \\

With the above quick review of asymptotic strong Feller property, we are now at the place to prove Theorem \ref{t:ModLogHarASF} and Corollary \ref{p:ModLogHarHaiMat}.
\Bp[Proof of Theorem \ref{t:ModLogHarASF}]
For any bounded differentiable function $f$, we choose some constant $\e>0$ small enough to make
$\e ||f||_\infty<1/2$. Applying the modified log-Harnack inequality \eqref{e:ModLogHar} to $2+2\e f$,  there exist some constants
$\alpha>0, \beta \geq 0,  C>0$, $\tilde C>0$ and some function $\delta(t) \ge 0$ with $\lim \limits_{t \rightarrow \infty} \delta(t)=0$ such that
\Bes
\begin{split}
P_t \log (2+2\e f)(y) \leq & \log P_t (2+2\e f)(x)+C|x-y|^\alpha \\
&\ +\delta (t)\tilde C |x-y|^\beta ||D\log (2+2\e f)||_\infty,
\end{split}
\Ees
which clearly implies
\Bes
\begin{split}
P_t \log (1+\e f)(y) \leq & \log \left(1+\e P_tf(x)\right)+C|x-y|^\alpha \\
&\ +2 \e \delta(t) \tilde C |x-y|^\beta ||Df||_\infty.
\end{split}
\Ees
Since $\e ||f||_\infty<1/2$, we have by Taylor expansion of the function $\log(1+x)$
\Be \label{e:Ptfy-Ptfx}
\e \big[P_t f(y)-P_t f(x)\big] \leq \e^2 ||f||^2_\infty+C |x-y|^{\alpha}+2 \e \delta(t) \tilde C |x-y|^\beta ||Df||_\infty.
\Ee
Dividing by $\e$ on the both side of the above inequality and exchanging $x$ and $y$, one has
\Be \label{e:PtEps}
|P_t f(y)-P_t f(x)| \leq \e ||f||^2_\infty+\frac{C|x-y|^{\alpha}}{\e}+2 \delta(t) \tilde C|x-y|^\beta ||Df||_\infty
\Ee
for any bounded differentiable $f$.

\vskip 1.2 mm
Next, we follow the idea in \cite[Proposition 3.12]{HM06}. For any $\gamma>0$, we define the metric
$d_{\gamma}(x,y)=1 \wedge \frac{1}{\gamma} |x-y|$ for any $x,y \in \mcl X$. It is clear that $||d_\gamma||_\infty:=\sup_{x, y \in \mcl X} d_\gamma(x, y) \leq 1$.
For any $f$ differentiable function with $||f||_{d_{\gamma}} \leq 1$ and $||f||_\infty \leq ||d_\gamma||_\infty$ (recall
Remark \ref{r:DuaMu1-Mu2}), by $||f||_{d_{\gamma}} \leq 1$ one has
$||Df||_\infty \leq \frac{1}{\gamma}$. \eqref{e:PtEps} implies
\Bes
\begin{split}
&\ \  |\int_{\mcl X} f(z) P_t(x,dz)-\int_{\mcl X} f(z) P_t(y,dz)| \leq \e+\frac{C|x-y|^{\alpha}}{\e}+\frac{2\delta(t) \tilde C|x-y|^\beta}{\gamma}.
\end{split}
\Ees
Since each bounded $f$ with $||f||_{d_\gamma}<\infty$ can be approximated by bounded differentiable function sequences, 
the above inequality and Lemma
\ref{l:DuaD} implies
 \Bes
\begin{split}
||P_n(x,\cdot)-P_n(y,\cdot)||_{d_{\gamma}} \leq \e+\frac{C|x-y|^{\alpha}}{\e}+\frac{2\delta(t) \tilde C|x-y|^\beta}{\gamma}
\end{split}
\Ees
Taking $\gamma=\sqrt{\delta(n)}$ and $\e=|x-y|^{\alpha/2}$, we have
\Be
||P_n(x,\cdot)-P_n(y,\cdot)||_{d_{\gamma_n}} \leq (1+C)|x-y|^{\frac \alpha 2}+2 \sqrt{\delta(n)} \tilde C |x-y|^{\beta}.
\Ee
which, by \eqref{e:ASFDef2}, immediately implies that $(P_t)_{t \geq 0}$
is asymptotically strong Feller at $x$.
\Ep
\ \ \ \
\Bp [Proof of Corollary \ref{p:ModLogHarHaiMat}]
From \eqref{e:PtEps}, taking $\e=|x-y|$ and letting $y \rightarrow x$, we immediately obtain
\eqref{e:HaiMatIne1}.
\Ep
\section{An example and Proof of Theorem \ref{t:SNSModLogHar}} \label{SExa}
In this section, we shall study 2D stochastic Navier-Stokes equation driven by highly degenerate but \emph{essentially elliptic} noises
as an example satisfying our modified log-Harnack inequality.
2D stochastic Navier-Stokes equation has been intensively studied in \cite{HM06}, \cite{Mat02}, \cite{EMS01}, \cite{FM95}, \cite{KS10} and the references therein. \\

Let us first give a quick introduction to the background
of 2D stochastic Navier-Stokes equations.
\subsection{2D stochastic Navier-Stokes systems} \label{s:SNSDeg}
Let $\T^2=(\R/2\pi)^2$ and let
\begin{align*}
& L^2_0(\T^2,\R^2)=\{x \in L^2(\T^2, \R^2); \int_{\T^2} x(\xi) d\xi=0\}, \\
& H=\{x \in L_0^2(\T^2, \R^2); \ div x=0\}.
\end{align*}
Moreover,
$$|\cdot| {\rm \ \ and} \ \ \Ll \cdot, \cdot\Rr {\rm \ denote \ the \ norm \ and \ the \ inner \ product \ of}\  H \ {\rm respectively}.$$
Let $\mcl P: L_0^2(\T^2, \R^2) \rightarrow H$ be the orthogonal
projection. Define the Stokes operator by
$$A =\mathcal{P}(-\Delta)$$
with $\Delta$ being the Laplacian on $L_0(\T^2,\R^2)$ and $D(A)=H^2(\T^2,{\R}^2) \cap H$. It is well known that
$\{e_k=\frac 1{2\pi}e^{i k \cdot x}: k \in \Z^2 \setminus \{0\}, x\in \T^2\}$ is an orthonormal basis of $L_0^2(\T^2, \R)$ and that
$\Delta$ is self-adjoint with the spectrum $\{-|k|^2: k \in \Z^2 \setminus \{0\}\}$. It is also clear that $\Delta e_k=-|k|^2 e_k$. For any
real number $\alpha$, one can define
$(-\Delta)^{\alpha}$ by the spectral decomposition as
$$(-\Delta)^{\alpha}=\sum_{k \in \Z^2 \setminus \{0\}} |k|^{2 \alpha} e_k \otimes e_k.$$
We can define the $\alpha$ order Stokes operator by
$$A^{\alpha}=\mcl P (-\Delta)^{\alpha}$$
with the domain defined by \eqref{e:AalpDom} below.

Under the orthonormal basis $\{e_k\}_{k \in \Z^2 \setminus \{0\}}$, $H$ can also be defined by
$$H=\left \{x=\sum_{k \in \Z^2 \setminus \{0\}} x_k e_k: x_k \in \R^2, k \cdot x_k=0, \sum_{k \in \Z^2 \setminus \{0\}}|x_k|^2<\infty \right\}.$$ Furthermore, we define
\Be \label{e:AalpDom}
D(A^{\alpha})=\left \{x=\sum_{k \in \Z^2\setminus \{0\}} x_k e_k: x_k \in \R^2, k \cdot x_k=0, \sum_{k \in \Z^2 \setminus \{0\}} |k|^{4 \alpha}|x_k|^2<\infty \right\}.
\Ee
It is clear that if $\alpha>0$ we have the following Poincare inequality
$$|x| \leq |A^\alpha x|$$
for any $x \in D(A^\alpha)$.  \\

We shall study the following highly degenerate 2D stochastic Navier-Stokes type equation
\begin{equation} \label{e:XEqn}
\begin{cases}
dX(t) + [\nu A X(t) + B(X(t),X(t))] dt = QdW_t, \\
X(0) =x,
\end{cases}
\end{equation}
where
\begin{itemize}
\item $\nu>0$ is the viscosity constant.
\item The nonlinear term $B$ is defined by
\begin{equation*} \label{e:NonlinearB}
 B(u,v)=\mathcal{P}[(u \cdot\nabla)v], \ B(u)=B(u,u) \ \ \ \forall \ u, v \in H^1(\T^d,{\R}^d) \cap H.
\end{equation*}
\item $W_t$ is the cylindrical Brownian motion on $H$ and $Q$ satisfies the highly degenerate condition as in Assumption \ref{a:Q}.
\end{itemize}
\ \

Given a $N \in \N$, define a project map $\pi_N: H \rightarrow H$ as follows:
for any $x \in H$ with $x=\sum_{k \in \Z^2 \setminus \{0\}} x_k e_k$, define
$$\pi_N x=\sum_{|k| \leq N} x_k e_k.$$
We split the space $H$ into the low and high frequency parts as
$$H=\pi_N H+(Id-\pi_N)H$$
where $Id$ is the identity map. For the generic $N \in \N$ we write $$H^\pl:=\pi_N H, \ \ H^\ph:=(Id-\pi_N) H.$$
For any $x \in H$, $x^\pl:=\pi_N x$ and $x^\ph:=(Id-\pi_N) x$.
It is clear that for any $\alpha>0$ one has
\Be \label{e:HLInq}
|A^{\alpha} x^{\pl}| \leq N^{2\alpha} |x^\pl|, \ \ \  |A^{\alpha} x^\ph| \geq N^{2 \alpha} |x^\ph|.
\Ee
\begin{assumption} [Highly degenerate but essentially elliptic noises assumption] \label{a:Q}
There exists some fixed $N_0 \in \N$ such that $Ran(Q)=H^\pl:=\pi_{N_0} H$ and
$Qx=0$ for any $x \in H^\ph$.
\end{assumption}

\BR
From this assumption, we clearly have $tr(QQ^{*})<\infty$ and that the operator $Q: H^\pl \rightarrow H^\pl$ is invertible, i.e. there exists some $C_0>0$ such that
\Be \label{e:Q-1c0}
|Q^{-1} x| \leq C_0 |x|
\Ee
for any $x \in H^\pl$. In our proof, we shall choose some large (but fixed) $N_0$ to make the noises
$Q dW_t$ have \emph{essential ellipticity} effect (see section 5.4 of \cite{HM06}).
\ER

Let us now write Theorem \ref{t:SNSModLogHar} in an exact form as the following Theorem \ref{t:ModLogHI},
which will be proven in next section.
\BT \label{t:ModLogHI}
There exist some $C=C(|x|,|y|)>0$ and $\tilde C=\tilde C(|x|,|y|)>0$ such that as $\nu N^2_0 > \frac 12 tr(QQ^{*})$ and $\nu>\max\{tr(QQ^{*}),C_2\}$ with $C_2>0$ defined in \eqref{e:X12Y12Z12} below, for any bounded differentiable function $f \geq 1$ we have
\Bes
\begin{split}
P_t \log f(y) \leq & \log P_tf(x)+C(|x-y|^2+|x-y|^4)\\
& \ + e^{-(\nu N_0^2-\frac 12 tr(QQ^*))t} \tilde C |x-y| ||D \log f||_\infty.
\end{split}
\Ees
The exact values of $C$ and $\tilde C$ can be easily figured out from the proof. Moreover,
\Bes
|DP_tf(x)| \leq ||f||^2+C+e^{-(\nu N_0^2-\frac 12 tr(QQ^*))t} \tilde C ||D \log f||_\infty.
\Ees
\ET

\ \ \ \ \ \
\subsection{Proof of Theorem \ref{t:ModLogHI}}
We shall apply asymptotic coupling method in the spirit of the idea in \cite[Proposition 4.11]{HM06}. For the more application of this method, we refer to \cite{EMS01}, \cite{Hai02}, \cite{Mat02} and \cite{HMS09}. \\

Our application of the asymptotic coupling method is sketched as follows. Give any $v \in L_{loc}^2([0,\infty);H)$ adapted to $\mcl F_t:=\sigma(W_s; 0 \leq s \leq t)$, define
\Be \label{e:Gir1}
\tilde W_t=W_t+\int_0^t v_s ds,
\Ee
by Girsanov theorem, we have a new probability measure $\tilde P$ under which
$\tilde W$ is a Brownian motion. This probability $\tilde P$ is uniquely determined by
\Be \label{e:Gir2}
\frac{d\tilde P}{dP}|_{\mcl F_t}=\exp\left\{-\int_0^t v_s dW_s-\frac12 \int_0^t |v_s|^2 ds\right\}.
\Ee

 Consider the SPDE
\begin{equation} \label{e:YEqn}
\begin{cases}
dY(t)+[A Y(t)+B(Y(t))]dt=Qd \tilde W_t, \\
Y(0)=y.
\end{cases}
\end{equation}
Denote $Z(t)=Y(t)-X(t)$, it is easy to see that
\begin{equation} \label{e:ZEqn}
\begin{cases}
\p_t Z(t)+A Z(t)+B(Z(t))+\tilde B(Z(t), X(t))=Q v_t, \\
Z(0)=z.
\end{cases}
\end{equation}
where $\tilde B(x,y)=B(x,y)+B(y,x)$ and $z=y-x$.

Eq. \eqref{e:ZEqn} can be split into two pieces, i.e. low frequency
and high frequency dynamics as the following
\begin{equation} \label{e:ZEqnLow}
\p_t Z^\pl (t)+A Z^\pl (t)+B^\pl (Z(t))+\tilde B^\pl (Z(t), X(t))=Q v_t
\end{equation}
with $Z^\pl (0)=z^\pl$, and
\begin{equation} \label{e:ZEqnHig}
\p_t Z^\ph (t)+A Z^\ph (t)+B^\ph (Z(t))+\tilde B^\ph (Z(t), X(t))=0
\end{equation}
with $Z^\ph (0)=z^\ph$. \\

Let us choose the $v$ in the following way. First of all, let \Be \label{e:Zlt}
\begin{cases}
Z^\pl(t)=(1-t)z^\pl \ \ \ \ 0 \leq t \leq 1,  \\
Z^\pl(t)=0 \ \ \ t>1.
\end{cases}
\Ee
Plugging this $Z^\pl(t)$ into \eqref{e:ZEqnHig}, we obtain the following PDE with unknown
$Z^\ph$
\Bes
\begin{cases}
\p_t Z^\ph (t)+A Z^\ph (t)+B^\ph (Z^\pl(t)+Z^\ph(t))+\tilde B^\ph (Z^\pl(t)+Z^\ph(t), X(t))=0, \\
Z^\ph(0)=z^\ph,
\end{cases}
\Ees
which has a unique solution by the same method as in \cite[Theorem 3.2]{Te95}.

Now $Z(t)=Z^\pl(t)+Z^\ph(t)$ is known. From Eq. \eqref{e:ZEqnLow}, we have
\Be \label{e:vCon}
v_t=
\begin{cases}
Q^{-1}[-z^\pl+(1-t)A z^{\pl}+B^\pl(Z(t))+\tilde B^\pl(Z(t),X(t))] \ \ \ 0 \leq t <1,\\
Q^{-1} [B^\pl(Z^\ph(t))+\tilde B^\pl(Z^\ph(t),X(t))] \ \ \ \ t \geq 1.
\end{cases}
\Ee
By the relation $Z=Y-X$, we also have
\Be \label{e:vCon1}
v_t=
\begin{cases}
Q^{-1}[-z^\pl+(1-t)A z^{\pl}-B^\pl(Z(t))+\tilde B^\pl(Z(t),Y(t))] \ \ \ 0 \leq t <1,\\
Q^{-1} [-B^\pl(Z^\ph(t))+\tilde B^\pl(Z^\ph(t),Y(t))] \ \ \ \ t \geq 1.
\end{cases}
\Ee
\ \\

To prove Theorem \ref{t:ModLogHI}, we need the following auxiliary lemmas which will be proven in the next section.
\BL \label{l:EntIne}
Let $f \geq 0$ with $\E f>0$. Then for any measurable function $g$, we have
\Bes
\E [f g] \leq \E f \log \E e^g+\E[f \log f]-\E f\log \E f
\Ees
\EL

\BL \label{l:XEst}
As $\nu>2tr(QQ^{*})$, we have \Bes \E_P \exp\left(|X(t)|^2+\nu
\int_0^t |A^{\frac 12}X(s)|^2 ds\right) \leq e^{|x|^2+tr(QQ^{*})t}.
\Ees
\Bes \E_{\tilde P} \exp\left(|Y(t)|^2+\nu
\int_0^t |A^{\frac 12}Y(s)|^2 ds\right) \leq e^{|y|^2+tr(QQ^{*})t}.
\Ees
\EL

\BL \label{l:B}
We have
\Be \label{e:DivFreRes}
\Ll x,B(y,x)\Rr=0, \ \ \ \Ll x, B(y,z)\Rr=-\Ll z,B(y,x)\Rr
\Ee
for all $x,y\in D(A^{\frac12})$ and $z \in D(A^{\frac 12})$. Furthermore, we have
\Be  \label{e:ZInf}
|\Ll x,B(y,z) \Rr| \leq C_1|x||y||A^{\frac 32}z|
\Ee
for all $x, y \in H$ and $z \in D(A^{\frac 32})$, and
\Be \label{e:X12Y12Z12}
|\Ll x,B(y,z) \Rr| \leq C_2 |x|^{1/2}|A^{\frac 12}x|^{1/2} |y|^{1/2}|A^{\frac 12}y|^{1/2} |A^{\frac 12}z|.
\Ee
for all $x,y,z \in D(A^{\frac 12})$. The constants $C_1, C_2$ both only depend on the space dimension.
\EL

\BL \label{l:BLEst}
For any $x, y \in D(A^{\frac 12})$, we have
\Be \label{e:BLEst}
|B^{\pl}(x,y)| \leq C_1 N_0^3 |x||y|
\Ee
where $C_1$ is the same as that in Lemma \ref{l:B}.
\EL

\BL
\label{l:HigParDec}
Let $v$ be chosen as in \eqref{e:vCon}. For any $p \geq 1$, if $\nu>\max\{C_2 \sqrt{p/2},2 tr(QQ^{*})\}$ we have
\Be \label{e:Zh0t1Est}
\begin{split}
& \E_P \sup_{0 \leq t \leq 1}|Z^\ph (t)|^{2p} \leq K_p e^{|x|^2} |x-y|^{2p}, \\
& \E_{\tilde P} \sup_{0 \leq t \leq 1}|Z^\ph (t)|^{2p} \leq K_p e^{|y|^2} |x-y|^{2p};
\end{split}
\Ee
and for $t>1$
\Be \label{e:Zht>1Est}
\begin{split}
& \E_P |Z^\ph (t)|^{2p} \leq \exp \LBr -(2 \nu p N_0^2-tr(QQ^{*})) t \RBr K_p e^{2|x|^2+2\nu p N_0^2} |x-y|^{2p}, \\
& \E_{\tilde P} |Z^\ph (t)|^{2p} \leq \exp \LBr -(2 \nu p N_0^2-tr(QQ^{*})) t \RBr K_p e^{2|y|^2+2\nu p N_0^2} |x-y|^{2p},
\end{split}
\Ee
where
\Bes
\begin{split}
K_p=& 2^{p-1} \exp\left\{C_1 p N_0^2 \big(|x-y|^2+|x-y|\big)+\frac{C_1 p N_0^3}2+tr(QQ^{*})\right\} \\
 & \times \left[\left(1+C_1 N^3_0+\frac {\nu N^2_0} 4\right)^p+p! \left(\frac{C^2_2}{4\nu}+\frac{C_1 N_0^3}2\right)^p \left(\frac{C^2_2 p}{4\nu}\right)^{-p}\right]
\end{split}
\Ees
and $C_1, C_2$ are the same as in Lemma \ref{l:B}.
\EL
\ \ \ \

\Bp [Proof of Theorem \ref{t:ModLogHI}] The second inequality in the theorem follows from Corollary \ref{p:ModLogHarHaiMat}
immediately. Let us prove the first inequality in the following two steps. \\

\emph{Step 1}.  Recalling Eq. \eqref{e:YEqn}, for any bounded differentiable function $f \geq 1$ one has
\Bes
P_t \log f(y)=\E_{\tilde P} \log f(Y(t))=\E_{\tilde P} \log f(Y^\pl(t)+Y^\ph(t)).
\Ees
By \eqref{e:Zlt} we have $Z^\pl(t)=0$, i.e. $X^\pl(t)=Y^\pl(t)$ for all $t \geq 1$, hence
\Bes
P_t \log f(y)=\E_{\tilde P} \Lbr \log f(X^\pl(t)+Y^\ph(t)) \Rbr, \ \ \ t \geq 1.
\Ees
Writing $z=x-y$, by Lemmas \ref{l:HigParDec} and \ref{l:EntIne} we have
\Be \label{e:MaiPro1}
\begin{split}
P_t \log f(y) &=\E_{\tilde P} \left[\log f(X^\pl(t)+Y^\ph(t))-\log f(X^\pl(t)+X^\ph(t))\right]+\E_{\tilde P} \log f(X(t)) \\
& \leq ||D \log f||_\infty \left\{\E_{\tilde P}|Z^\ph(t)|^2\right\}^{\frac 12}+\E_P[\frac{d \tilde P}{dP} \log f(X(t))] \\
& \leq \exp \LBr -(\nu N_0^2-\frac 12 tr(QQ^{*})) t+|y|^2+\nu  N_0^2\RBr
\sqrt {K_1}  |z| ||D \log f||_\infty \\
& \ \ +\E_P\left[\frac{d \tilde P}{dP} \log \frac{d \tilde P}{dP}\right]+\log \E_P f(X(t))
\end{split}
\Ee
where $K_1$ is defined in Lemma \ref{l:HigParDec}.
For the entropy term above, by \eqref{e:Gir1} and \eqref{e:Gir2},
\Be \label{e:MaiPro2}
\begin{split}
\E_P\left[\frac{d \tilde P}{dP}  \log \frac{d \tilde P}{dP}\right]&=\E_{\tilde P}\left[-\int_0^t v_s dW_s-\frac12 \int_0^t |v_s|^2ds\right] \\
&=\E_{\tilde P}\left[-\int_0^t v_s d\tilde W_s+\frac12 \int_0^t |v_s|^2ds\right]=\frac 12\E_{\tilde P}\int_0^t |v_s|^2ds.
\end{split}
\Ee
We claim
\Be \label{e:vsEst}
\E_{\tilde P}\int_0^t |v_s|^2ds \leq (L_1+L_3)|z|^4+(L_2+L_4)|z|^2.
\Ee
where $L_1 \cdots, L_4$ are defined in Step 2 below.
From \eqref{e:MaiPro1}, \eqref{e:MaiPro2} and \eqref{e:vsEst}, we have
\Bes
\begin{split}
P_t \log f(y) & \leq \log P_t f(x)+(L_1+L_3)|z|^4+(L_2+L_4)|z|^2 \\
&+\exp \LBr -(\nu N_0^2-\frac 12 tr(QQ^{*})) t+|y|^2+\nu p N_0^2\RBr
\sqrt {K_1}  |z| ||D \log f||_\infty.
\end{split}
\Ees
By the definitions of $K_1, L_1, \cdots L_4$ and recalling $z=y-x$, we conclude the proof up to proving \eqref{e:vsEst}, and can easily figure out the exactly values of $C$ and $\tilde C$ in the theorem. \\

\emph{Step 2}. Let us prove \eqref{e:vsEst}. We first consider $\E_{\tilde P}\int_0^1 |v_s|^2 ds$.
 By \eqref{e:Q-1c0}, \eqref{e:vCon1} and $|Az^\pl| \leq N_0^2 |z|$, one has
\Bes
\E_{\tilde P}\int_0^1 |v_s|^2ds \leq 3 C^2_0 \left(4 N^4_0|z|^2+\E_{\tilde P}\int_0^1|B^\pl(Z(s))|^2+|\tilde B^\pl(Z(s),Y(s))|^2 ds \right)
\Ees
By Lemma \ref{l:BLEst}, \eqref{e:Zh0t1Est} and $|Z^\pl(t)| \leq |z|$ for $0 \leq t \leq 1$, we have
\Bes
\begin{split}
& \ \E_{\tilde P}\int_0^1 |B^\pl(Z(s))|^2ds \leq C^2_1 N^6_0 \int_0^1 \E_{\tilde P} |Z(s)|^4 ds \\
&=8 C^2_1 N^6_0 \int_0^1 \E_{\tilde P} \left(|Z^\pl(s)|^4+|Z^\ph(s)|^4\right) ds  \\
& \leq 8C^2_1 N^6_0 \big(1+K_2 e^{|y|^2}\big)|z|^4.
\end{split}
\Ees
Moreover, by Lemmas \ref{l:BLEst} and \ref{l:XEst}, \eqref{e:Zh0t1Est} and a similar argument as the above, one has
\Bes
\begin{split}
&\int_0^1 |\tilde B^\pl(Z(s),Y(s))|^2 ds \leq 4 C^2_1 N^6_0 \int_0^1 \E_{\tilde P} |Z(s)|^2|Y(s)|^2ds \\
 & \leq 4 C^2_1 N^6_0 \int_0^1 \left(\E_{\tilde P} Z^4(s)\right)^{\frac 12} \left(\E_{\tilde P} |Y(s)|^4\right)^{\frac 12}ds \\
 &\leq 4 C^2_1 N^6_0 \int_0^1 \left(\E_{\tilde P} Z^4(s)\right)^{\frac 12} \left(2\E_{\tilde P} e^{|Y(s)|^2}\right)^{\frac 12}ds\\
& \leq 4 \sqrt 2 C^2_1 N^6_0 \sqrt{1+K_2 e^{|y|^2}} \exp\LBr \frac 12|y|^2+\frac 12tr(QQ^{*}) \RBr |z|^2.
\end{split}
\Ees
Collecting all the above, we have
\Be \label{e:EP01vs}
\E_{\tilde P}\int_0^1 |v_s|^2ds  \leq L_1|z|^4+L_2|z|^2.
\Ee
where
$$L_1=24C^2_0C^2_1 N^6_0 \big(1+K_2 e^{|y|^2}\big),$$
$$L_2=3C^2_0 \lbr 4 N^4_0+ 4 \sqrt 2 C^2_1 N^6_0 \sqrt{1+K_2 e^{|y|^2}} \ e^{[|y|^2+tr(QQ^{*})]/2}\rbr.$$

Now let us estimate $\E_{\tilde P} \int_1^t |v_s|^2 ds$, by \eqref{e:vCon1} one has
\Bes
\E_{\tilde P}\int_1^t |v_s|^2ds \leq 2 C^2_0 \E_{\tilde P}\int_1^t \left(|B^\pl(Z^\ph (s))|^2+|\tilde B^\pl(Z^\ph (s),Y(s))|^2 \right)ds.
\Ees
By a similar argument as for proving \eqref{e:EP01vs} and thanks to
\eqref{e:Zht>1Est}, when $\nu N^2_0>\frac 12 tr(QQ^{*})$ one has
\Bes
\E_{\tilde P} \int_1^t |B^\pl(Z^\ph (s))|^2 ds \leq C^2_1 N^6_0 \int_1^t \E_{\tilde P} |Z^\ph (s)|^4 ds
\leq \frac{C^2_1 N^6_0 e^{2|y|^2+4 \nu N^2_0}K_2}{4 \nu N^2_0-tr(QQ^{*})} |z|^4,
\Ees
and
\Bes
\begin{split}
& \ \ \ \ \E_{\tilde P}\int_1^t |\tilde B^\pl(Z^\ph (s),Y(s))|^2ds \leq 2C^2_1 N^6_0 \int_1^t \E_{\tilde P}|Z^\ph(s)|^2 |Y(s)|^2 ds \\
& \leq 2C^2_1 N^6_0 \int_1^t \lbr \E_{\tilde P}|Z^\ph(s)|^4 \rbr^{1/2} \lbr \E_{\tilde P}|Y(s)|^4 \rbr^{1/2} ds \leq \frac{2 C_1^2 N_0^6\sqrt{2K_2} e^{\frac 32|y|^2+2 \nu N^2_0}}{2 \nu N^2_0-tr(QQ^{*})} |z|^2
\end{split}
\Ees
Therefore,
\Be \label{e:EP1tIntVsEst}
\E_{\tilde P}\int_1^t |v_s|^2ds \leq L_3 |z|^4+L_4 |z|^2
\Ee
where
$$L_3=\frac{2C^2_0C^2_1 N^6_0 e^{2|y|^2+4 \nu N^2_0}K_2}{4 \nu N^2_0-tr(QQ^{*})}, \ \ L_4=\frac{4C^2_0 C_1^2 N_0^6\sqrt{2K_2} e^{\frac 32|y|^2+2 \nu N^2_0}}{2 \nu N^2_0-tr(QQ^{*})}.$$
\Ep

\section{Proof of auxiliary lemmas in section \ref{SExa}}
Some of the first four lemmas are well known. Since their proofs are short, it is very convenient to repeat them here.
\Bp [Proof of Lemma \ref{l:EntIne}]
Since all the expectations restricted on the
set $\{x: f(x)=0\}$ are zero, without loss of generality, we can assume that $f>0 \ a.e.$.
We can also simply assume that $\E f=1$, otherwise one can replace
$f$ by $\frac{f}{\E f}$. We have
\Bes
\begin{split}
\E[f g] \leq \E[f \log e^g]& =\E[f \log \frac{e^g}{f}]+\E[f \log f] \\
& \leq \log \E [f \frac{e^g}{f}]+\E[f \log f]
=\E e^g+\E[f \log f].
\end{split}
\Ees
\Ep

\Bp[Proof of Lemma \ref{l:XEst}] The proofs for the two claims are the same, so we only prove the first one.
 By It$\hat o$ formula, we have \Bes
 \begin{split}
 & \ \ |X(t)|^2+\nu
\int_0^t |A^{\frac 12} X(s)|^2 ds \\
&=|x|^2+tr(QQ^{*})t+2 \int_0^t \Ll
X(s),QdW_s\Rr-\nu
\int_0^t |A^{\frac 12} X(s)|^2 ds.
\end{split}
\Ees
By $|x| \leq |A^{\frac 12} x|$ and $|Q^{*} x|^2 \leq tr(QQ^{*}) |x|^2$, we have
\Bes
\begin{split}
& \ \ \E_P \exp\left\{2\int_0^t \Ll
X(s),QdW_s\Rr-\nu \int_0^t |A^{\frac 12} X(s)|^2 ds \right\} \\
& \leq \E_P \exp\left\{2 \int_0^t \Ll
X(s),QdW_s\Rr-2 \int_0^t |Q^{*} X(s)|^2ds+(2 tr(QQ^{*})-\nu)\int_0^t |X(s)|^2 ds \right\}.
\end{split}
\Ees
Since $\exp\left\{2 \int_0^t \Ll
X(s),QdW_s\Rr-2 \int_0^t |Q^{*} X(s)|^2ds\right\}$ is a martingale, as $\nu>2tr(QQ^{*})$,
one has
\Bes
\E_P \exp \left\{|X(t)|^2+\nu \int_0^t |A^{\frac 12} X(s)|^2 ds \right\} \leq e^{|x|^2+tr(QQ^{*})t}.
\Ees
\Ep

\Bp [Proof of Lemma \ref{l:B}]
\eqref{e:DivFreRes} is classical, one can, for instance, refers to \cite{Te95}.
We clearly have
$$|\Ll x,B(y,z)\Rr| \leq |x||y|||\nabla z||_\infty \leq C_1 |x||y||A^{\frac 32}z|$$
since $z=\sum \limits_{k \in \Z^{2}\setminus \{0\}} z_k e^{i 2\pi k \cdot x}$ and hence
\Bes
\begin{split}
||\nabla z||_\infty  \leq  C \sum_{k \in \Z^2 \setminus \{0\}} |k| |z_k| \leq C \sqrt{\sum_{k \in \Z^2 \setminus \{0\}}|k|^{-4}}\sqrt{\sum_{k \in \Z^2 \setminus \{0\}} |k|^6 |z_k|^2}
\end{split}
\Ees
As for \eqref{e:X12Y12Z12}, by H$\ddot{o}$lder's inequality, the classical Sobolev embedding
$||y||_{L^4} \leq C|A^{\frac 14} y|$ (\cite[Theorem 6.16 and Remark 6.17]{Hai09Note}) and the easy interpolation $|A^{\frac 14} x| \leq |x|^{1/2}|A^{\frac 12}x|^{1/2}$, we have
\Bes |\Ll x,B(y,z)\Rr|  \leq ||x||_{L^4} ||y||_{L^4} |A^{\frac 12}z|
\leq C_2|x|^{1/2}|A^{\frac 12}x|^{1/2} |y|^{1/2}|A^{\frac 12}y|^{1/2} |A^{\frac 12}z|.
\Ees
\Ep

\Bp [Proof of Lemma \ref{l:BLEst}]
For any $z \in H$ \emph{with} $z^\ph=0$, it is clear from the first inequality of \eqref{e:HLInq} that
$$|A^{\frac 32} z| \leq N_0^3 |z|.$$
By \eqref{e:DivFreRes} and \eqref{e:ZInf}, one has
\Bes
|\Ll z, B^\pl(x,y) \Rr|=|\Ll y, B(x,z) \Rr|
 \leq C_1 |A^{\frac 32} z||y||x| \leq C_1 N_0^3|y||x||z|
 \Ees
\Ep

\Bp [Proof of Lemma \ref{l:HigParDec}]
 The proof of the claims for $P$ and $\tilde P$ are the same, so we only show that for the former. We split the proof into the following three steps. \\

\emph{Step 1.} By \eqref{e:ZEqnHig} we have
\Be \label{e:EneZHig}
\p_t |Z^\ph(t)|^2+2 \nu |A^{\frac 12} Z^\ph(t)|^2=-2\Ll Z^\ph(t),B^\ph(Z(t))\Rr-2\Ll
Z^\ph(t),\tilde B^\ph (Z(t), X(t))\Rr. \Ee
By \eqref{e:DivFreRes}, \eqref{e:ZInf} and \eqref{e:HLInq} we have
\Be \label{e:ZhBZ}
\begin{split}
\Ll Z^\ph,B^\ph(Z)\Rr&=\Ll Z^\ph,B(Z^\ph)+\tilde B(Z^\ph,Z^\pl)+B(Z^\pl)\Rr \\
&=\Ll Z^\ph, B(Z^\ph,Z^\pl)\Rr+\Ll Z^\ph,B(Z^\pl)\Rr
\\
& \leq C_1 N_0^3 (|Z^\pl||Z^\ph|^2+|Z^\pl|^2|Z^\ph|).
\end{split}
\Ee
As for the second term on the r.h.s. of \eqref{e:EneZHig}, we have by \eqref{e:X12Y12Z12}
\Be \label{e:ZhBZX}
\begin{split}
|\Ll Z^\ph, B (Z, X)\Rr|& \leq C_2|Z^\ph|^{1/2} |A^{\frac 12}Z^\ph|^{1/2} |Z|^{1/2} |A^{\frac 12} Z|^{1/2} |A^{\frac
12} X|,
\end{split}
\Ee
and have by \eqref{e:DivFreRes}, \eqref{e:ZInf} and \eqref{e:HLInq}
\Be \label{e:ZhBXZ}
\begin{split}
|\Ll Z^\ph, B(X, Z)\Rr| &=|\Ll Z^\ph, B(X,Z^\pl)\Rr| \leq C_1 N^3_0 |Z^\ph||X||Z^\pl|.
\end{split}
\Ee

\emph{Step 2.} Let us now estimate $\E \sup_{0 \leq t\leq 1} |Z^\ph(t)|^{2p}$. By
\eqref{e:ZhBZ}, Cauchy inequality and $|Z^\pl(t)| \leq |z|$ for $0 \leq t \leq 1$ (see \eqref{e:Zlt}) one has
\Bes
\begin{split}
|\Ll Z^\ph, B^\ph (Z)\Rr| &\leq C_1N^3_0 |z| |Z^\ph|^2+\frac 12 C_1 N^3_0 |z|^2|Z^\ph|^2+\frac 12 C_1 N^3_0|z|^2 \\
&\leq C_1N^3_0 \big(|z|^2+|z|\big)|Z^\ph|^2+C_1 N^3_0|z|^2.
\end{split}
\Ees
Applying Young's inequality (two times), \eqref{e:ZhBZX} and $|A^{\frac 12} z^{\pl}| \leq N_0 |z|$, we have
\Bes
\begin{split}
|\Ll Z^\ph, B (Z, X)\Rr| & \leq \frac{\nu} 2 |A^{\frac 12}Z^\ph||A^{\frac 12} Z|+\frac{C^2_2}{4\nu}|A^{\frac 12} X|^2|Z^\ph|^2+\frac{C^2_2}{4 \nu}|A^{\frac 12} X|^2|Z|^2 \\
& \leq
\frac{3\nu}{4}|A^{\frac 12} Z^\ph|^2+\frac \nu 4|A^{\frac 12} Z^\pl|^2+\frac{C^2_2}{4\nu}|A^{\frac 12}X|^2 |Z^\ph|^2+\frac{C^2_2}{4\nu}|A^{\frac 12}X|^2|Z^\pl|^2 \\
&\leq \frac{3\nu}{4}|A^{\frac 12} Z^\ph|^2+\frac {\nu N^2_0} 4 |z|^2+\frac{C^2_2}{4\nu}|A^{\frac 12}X|^2 |Z^\ph|^2+\frac{C^2_2}{4\nu}|A^{\frac 12}X|^2|z|^2.
\end{split}
\Ees
By \eqref{e:ZhBXZ}, $|X| \leq |A^{\frac 12} X|$ and $|Z^{\pl}(t)| \leq |z^\pl|$ for $0 \leq t\leq 1$, we get
\Bes
\begin{split}
|\Ll Z^\ph, B(X, Z)\Rr|\leq \frac{C_1 N_0^3}2 |Z^\ph|^2+\frac{C_1 N_0^3}2 |A^{\frac 12} X|^2|z|^2.
\end{split}
\Ees
The above three inequalities and \eqref{e:EneZHig} imply
\Bes
\begin{split}
\p_t |Z^\ph(t)|^2 &\leq \left[C_1 N_0^2 \big(|z|^2+|z|\big)+\frac{C^2_2}{4\nu}|A^{\frac 12}X|^2+ \frac{C_1 N_0^3}2\right] |Z^\ph(t)|^2 \\
&\ \ +\left[C_1 N^3_0 |z|^2+\frac {\nu N^2_0} 4 |z|^2+\frac{C^2_2}{4\nu}|A^{\frac 12}X|^2|z|^2+\frac{C_1 N_0^3}2 |A^{\frac 12} X|^2|z|^2\right],
\end{split}
\Ees
thus  for any $0 \leq t \leq 1$
\Bes
\begin{split}
|Z^\ph(t)|^2 &\leq \exp\left\{C_1 N_0^2 \big(|z|^2+|z|\big)+\frac{C_1 N_0^3}2+\frac{C^2_2}{4\nu}\int_0^1|A^{\frac 12}X(s)|^2 ds \right\} \\
& \times \left[|z^\ph|^2+C_1 N^3_0 |z|^2+\frac {\nu N^2_0} 4 |z|^2+\left(\frac{C^2_2}{4\nu}+\frac{C_1 N_0^3}2\right)|z|^2 \int_0^1 |A^{\frac 12} X(s)|^2 ds\right].
\end{split}
\Ees
It follows that
\Bes
\begin{split}
\sup_{0 \leq t \leq 1}|Z^\ph(t)|^{2p} &\leq 2^{p-1} \exp\left\{C_1 p N_0^2 \big(|z|^2+|z|\big)+\frac{C_1 p N_0^3}2\right\}\exp\left\{\frac{C^2_2 p}{4\nu}\int_0^1|A^{\frac 12}X(s)|^2 ds \right\} \\
& \times \left[\left(1+C_1 N^3_0+\frac {\nu N^2_0} 4\right)^p+\left(\frac{C^2_2}{4\nu}+\frac{C_1 N_0^3}2\right)^p \left(\int_0^1 |A^{\frac 12} X(s)|^2 ds\right)^p\right] |z|^{2p}.
\end{split}
\Ees
By Lemma \ref{l:XEst}, as $\nu>\max\{C_2 \sqrt{p}/2, 2tr(QQ^{*})\}$,
\Bes
\E \exp\left\{\frac{C^2_2 p}{4\nu}\int_0^1|A^{\frac 12}X(s)|^2 ds \right\} \leq e^{|x|^2+tr(QQ^{*})}
\Ees
and
\Bes
\begin{split}
& \ \ \E \left[\exp\left\{\frac{C^2_2 p}{4\nu}\int_0^1|A^{\frac 12}X(s)|^2 ds \right\} \left(\int_0^1 |A^{\frac 12} X(s)|^2 ds\right)^p \right] \\
&\leq p! \left(\frac{C^2_2 p}{4\nu}\right)^{-p} \E \exp\left\{\frac{C^2_2 p}{4\nu}\int_0^1|A^{\frac 12}X(s)|^2 ds \right\}  \leq p! \left(\frac{C^2_2 p}{4\nu}\right)^{-p} e^{|x|^2+tr(QQ^{*})}.
\end{split}
\Ees
Collecting the above three inequalities, we get
\Bes
\begin{split}
\E\sup_{0 \leq t \leq 1}|Z^\ph(t)|^{2p} & \leq 2^{p-1} \exp\left\{C_1 p N_0^2 \big(|z|^2+|z|\big)+\frac{C_1 p N_0^3}2\right\} \\
& \times \left[\left(1+C_1 N^3_0+\frac {\nu N^2_0} 4\right)^p+p! \left(\frac{C^2_2}{4\nu}+\frac{C_1 N_0^3}2\right)^p \left(\frac{C^2_2 p}{4\nu}\right)^{-p}\right]  e^{|x|^2+tr(QQ^{*})} |z|^{2p}
\end{split}
\Ees
\ \ \

\emph{Step 3}. As $t>1$, $Z^\pl(t)=0$ by \eqref{e:Zlt}. From \eqref{e:EneZHig}-\eqref{e:ZhBXZ},
\Be \label{e:EneHigT>1}
\p_t |Z^\ph|^2+2 \nu |A^{\frac 12} Z^\ph|^2 \leq C_2|Z^\ph||A^{\frac 12} Z^\ph| |A^{\frac 12} X|,
\Ee
which, together with Young's inequality, implies
\Bes
\p_t |Z^\ph|^2+\nu |A^{\frac 12} Z^\ph|^2  \leq \frac{C^2_2}{4 \nu} |A^{\frac 12} X|^2 |Z^\ph|^2.
\Ees
By the second inequality of \eqref{e:HLInq} we further have
\Bes
\p_t |Z^\ph|^2+\nu N^2_0 |Z^\ph|^2  \leq \frac{C^2_2}{4 \nu} |A^{\frac 12} X|^2 |Z^\ph|^2.
\Ees
Therefore, for all $t>1$
\Bes
\begin{split}
\E |Z^\ph(t)|^{2p} & \leq \E \exp \LBr -\nu p N_0^2(t-1)+\frac{C^2_2 p}{4 \nu} \int_1^t |A^{\frac 12} X(s)|^2 ds \RBr |Z^\ph(1)|^{2p} \\
& \leq \left(\E \exp\left\{-2 \nu p N_0^2(t-1)+\frac{C^2_2 p}{2 \nu} \int_1^t |A^{\frac 12} X(s)|^2 ds\right\}\right)^{1/2} \left(\E |Z^\ph(1)|^{4p} \right)^{1/2} \\
& \leq \exp \LBr 2\nu p N_0^2+|x|^2-(2 \nu p N_0^2-tr(QQ^{*})) t \RBr  \left(\E |Z^\ph(1)|^{4p}\right)^{1/2}
\end{split}
\Ees
as $\nu>\max \{C_2 \sqrt{p/2},2 tr(QQ^{*})\}$, where the last inequality is due to Lemma \ref{l:XEst}. This, together with the last inequality in Step 2, immediately implies \eqref{e:Zht>1Est}
\Ep



\bibliographystyle{amsplain}


\end{document}